\documentclass[11pt]{article}
\usepackage{amssymb,amsfonts,amsmath,amsthm}
\usepackage{epsfig}
\usepackage{cite}
\newtheorem{theorem}{Theorem}[section]
\newtheorem{thm}[theorem]{Theorem}
\newtheorem{prop}[theorem]{Proposition}
\newtheorem{lem}[theorem]{Lemma}

\makeatletter \@addtoreset{equation}{section}

\begin{document}
\title{New extensions to the sumsets with polynomial restrictions}
\author{
  {\vspace{0.2cm}}
  {Yue Zhou}\\
 {\small School of Mathematical Science and Computing Technology}\\
{\small Central South University, Changsha 410075, P.R. China}\\
  {nkzhouyue@gmail.com}\\
 }

\date{February 14, 2012}
\maketitle

\begin{abstract}
By taking the leading and the second leading coefficients of the Morris identity,
we get new polynomial coefficients. These coefficients
lead to new results in the sumsets with polynomial restrictions
by the polynomial method of N. Alon.
\end{abstract}

\section{Introduction}

Let $A_1,\ldots,A_n$ be finite subsets of a field $F$ with
$0<k_1=|A_1|\leqslant \cdots \leqslant k_n=|A_n|$, where the characteristic of $F$
is infinite or a prime. We are concerning the lower bounds for
the following sumsets:
$$\{a_1+\cdots+a_n: a_1\in A_1,\ldots,a_n\in A_n, \ \ a_i 's \ \mbox{satisfy certain restrictions}\}.$$

In 1999, N. Alon got the following lemma which is called the Combinatorial
Nullstellensatz.
\begin{lem}\emph{(Alon \cite{Alon1})}
Let $A_1,\ldots,A_n$ be finite subsets of a field $F$
with $|A_i|>k_i$ for $i=1,\ldots,n$, where $k_1,\ldots,k_n$ are nonnegative integers. If the coefficient of the monomial $x_1^{k_1}\cdots x_n^{k_n}$ in $f(x_1,\ldots,x_n)\in F[x_1,\ldots,x_n]$ is nonzero and
$k_1+\cdots+k_n$ is the total degree of $f$, then there are $a_1\in A_1, \ldots, a_n\in A_n$ such that $f(a_1,\ldots,a_n)\neq 0$.
\end{lem}
The Combinatorial Nullstellensatz soon led to many
results \cite{Hou-Sun,Liu-Sun,Pan-Sun,Sun1,Sun07,Sun2,Sun-Yeh,Sun-Zhao} in the sumsets.
It also implies a polynomial method in finding the lower bounds for various
restricted sumsets, which is the following lemma.
\begin{lem}\label{l-polylem}
\emph{(Alon et al. \cite{ANR1,ANR2})}. Let $A_1, \ldots ,A_n$ be finite
nonempty subsets of a field $F$ with $|A_i| = k_i$ for $i = 1, \ldots , n.$ Let
$P(x_1, \ldots , x_n) \in F[x_1, \ldots , x_n] \backslash \{0\}$ and $\deg P \leqslant \sum_{i=1}^n (k_i - 1).$ If the
coefficient of the monomial $x_{1}^{k_1-1} \cdots x_{n}^{k_n - 1}$ in the polynomial
$$P(x_1, \ldots , x_n)(x_1 + \cdots + x_n)
^{\sum_{i=1}^n(k_i-1)-\deg P}$$
does not vanish, then we have
$$|\{a_1 + \cdots + a_n : a_i \in A_i, P(a_1, \ldots , a_n) \neq 0\}| \geqslant \sum_{i=1}^n (k_i - 1) - \deg P + 1.$$
\end{lem}

Throughout this paper, we need the following definitions and notations.

Let $k$, $m$ be nonnegative integers and let $n$ be a positive integer.
Denote $T$ by a nonempty subset of $\{1,\ldots,n\}$ and denote $F$ by a field with
characteristic $p$(infinite or a prime).
For $i, j = 1, \ldots , n$ with $i \neq j$,
let $A_i$ and $S_{ij}$ be subsets of $F$.
Define
\begin{itemize}
\item Condition (a): $|A_i|=k,|S_{ij}|\leqslant 2m$.

\item Condition (b): $|A_i|=k-n+i,|S_{ij}|< 2m$.
\end{itemize}

By using an equivalent form of Lemma \ref{l-polylem} (Proposition 2.1 of \cite{Hou-Sun}), Q. H. Hou and Z. W. Sun got the following theorem.
\begin{thm}\emph{(Hou and Sun \cite{Hou-Sun})}\label{thm-HouSun}
For $i, j = 1, \ldots , n$ with $i \neq j$, let $A_i$ and $S_{ij}$ satisfy
the \emph{Condition (a)}.
If $p>\max \{mn, (k-1)n-mn(n-1)\}$, then for the set
\begin{align}\label{e-setC}
 C = \{a_1 + \cdots + a_n: a_1 \in A_1, \cdots , a_n \in A_n, a_i - a_j \notin S_{ij} \ \ \mbox{if} \ \ i \neq j\},
\end{align}
we have
\begin{equation*}
|C| \geqslant (k + m - mn - 1)n + 1.
\end{equation*}
\end{thm}
In the original paper of Q. H. Hou and Z. W. Sun \cite{Hou-Sun},
they restrict $|S_{ij}|\leqslant m$, but the result is still right when $|S_{ij}|$ are extended to
$2m$, as pointed out by Z. W. Sun in \cite{Sun07}.

After that, Z. W. Sun and Y. N Yeh got a
theorem very closely to Theorem \ref{thm-HouSun}.
\begin{thm}\emph{(Sun and Yeh \cite{Sun-Yeh})}\label{thm-SunYeh}
For $i, j = 1, \ldots , n$ with $i \neq j$, let $A_i$ and $S_{ij}$ satisfy
the \emph{Condition (b)}.
If $p>\max \{mn, (k-1)n-mn(n-1)\}$, then we have
\begin{equation*}
|C| \geqslant (k + m - mn - 1)n + 1,
\end{equation*}
where $C$ is defined as in \eqref{e-setC}.
\end{thm}
In this paper,
we will show that Theorem \ref{thm-HouSun} and Theorem \ref{thm-SunYeh}
are actually equivalent by a relation between two kinds of polynomial coefficients.
By getting new polynomial coefficients we obtain some further results
related to Theorem \ref{thm-HouSun} and Theorem \ref{thm-SunYeh}.

By using the leading coefficient of the Morris identity \cite{morris82},
we get the following theorem.

\begin{thm}\label{thm1}
For $i, j = 1, \ldots , n$ with $i \neq j$, let $A_i$ and $S_{ij}$ satisfy
the \emph{Condition (a)}, and
let $\sum_{i\in T}\alpha_i\neq 0$ for $\alpha_i\in F$.
If $p>\max \{mn, (k-1)n-mn(n-1)-1\}$, then for the set
\begin{align}
 C = \{a_1 + \cdots + a_n: a_1 \in A_1, \cdots , a_n \in A_n, a_i - a_j \notin S_{ij}  \ \ \mbox{if} \ \ i \neq j
\ \  \mbox{and} \ \ \sum_{i\in T}\alpha_i a_i\neq 0\},
\end{align}
we have
\begin{equation*}
|C| \geqslant (k + m - mn - 1)n.
\end{equation*}
\end{thm}

By a relation between two kinds of polynomial coefficients
in Section 2, we
get the following theorem corresponds to Theorem \ref{thm1}.

\addtocounter{theorem}{-1}
\renewcommand{\thetheorem}{\thesection.\arabic{theorem}$'$}
\begin{thm}\label{thm1'}
For $i, j = 1, \ldots , n$ with $i \neq j$, let $A_i$ and $S_{ij}$ satisfy
the \emph{Condition (b)},
and let $C, \alpha_i$ and $p$ be as in Theorem \ref{thm1}.
We have
\begin{equation}
|C| \geqslant (k + m - mn - 1)n.
\end{equation}
\end{thm}

\renewcommand{\thetheorem}{\thesection.\arabic{theorem}}
By using the second leading coefficient of the Morris identity, we can get
the following theorems.
\begin{thm}\label{thm2}
For $i, j = 1, \ldots , n$ with $i \neq j$, let $A_i$ and $S_{ij}$ satisfy
the \emph{Condition (a)},
and let $\sum_{i\in T}\alpha_i\neq 0$ for $\alpha_i\in F$.
If $p>\max \{mn, (k-1)n-mn(n-1)-2\}$, then for the set
\begin{align}
 C = \{a_1 + \cdots + a_n: a_1 \in A_1, \cdots , a_n \in A_n, a_i - a_j \notin S_{ij}  \ \ \mbox{if} \ \ i \neq j
\ \  \mbox{and} \ \ \sum_{i\in T}\alpha_i a_i^2\neq 0\},
\end{align}
we have
\begin{equation*}
|C| \geqslant (k + m - mn - 1)n - 1.
\end{equation*}
\end{thm}

\addtocounter{theorem}{-1}
\renewcommand{\thetheorem}{\thesection.\arabic{theorem}$'$}
\begin{thm}\label{thm2'}
For $i, j = 1, \ldots , n$ with $i \neq j$, let $A_i$ and $S_{ij}$ satisfy
the \emph{Condition (b)},
and let $C, \alpha_i$ and $p$ be as in Theorem \ref{thm2}.
We have
\begin{equation*}
|C| \geqslant (k + m - mn - 1)n - 1.
\end{equation*}
\end{thm}

\renewcommand{\thetheorem}{\thesection.\arabic{theorem}}
\begin{thm}\label{thm3}
For $i, j = 1, \ldots , n$ with $i \neq j$, let $A_i$ and $S_{ij}$ satisfy
the \emph{Condition (a)}, and
let $\sum_{i,j\in T}\alpha_{ij}\neq 0$ for $\alpha_{ij}\in F$.
If $p>\max \{mn, (k-1)n-mn(n-1)-2\}$, then for the set
\begin{align}
 C = \{a_1 + \cdots + a_n: a_1 \in A_1, \cdots , a_n \in A_n, a_i - a_j \notin S_{ij}  \ \  \mbox{and} \ \ \sum_{i\in T}\alpha_{ij} a_ia_j\neq 0, \ \ \mbox{if} \ \ i \neq j\},
\end{align}
we have
\begin{equation}
|C| \geqslant (k + m - mn - 1)n - 1.
\end{equation}
\end{thm}

\addtocounter{theorem}{-1}
\renewcommand{\thetheorem}{\thesection.\arabic{theorem}$'$}
\begin{thm}\label{thm3'}
For $i, j = 1, \ldots , n$ with $i \neq j$, let $A_i$ and $S_{ij}$ satisfy
the \emph{Condition (b)},
and let $C, \alpha_{ij}$ and $p$ be as in Theorem \ref{thm3}.
We have
\begin{equation}
|C| \geqslant (k + m - mn - 1)n - 1.
\end{equation}
\end{thm}
\renewcommand{\thetheorem}{\thesection.\arabic{theorem}}

This paper is organized as follows. In Section 2 we give
a relation between two kinds of polynomial coefficients.
In Section 3 we obtain new polynomial
coefficients, which are the keys
to the polynomial method by N. Alon.
In Section 4 we prove Theorem \ref{thm1} with the polynomial method,
and the proof of other theorems are routine according to
the proof of Theorem \ref{thm1}.

\section{A relation between two kinds of polynomial coefficients}

Let
\begin{align}\label{e-H'}
H'_m(x)=\prod_{1\leqslant i < j \leqslant n}(x_i-x_j)^{2m}\cdot L(x_1,\ldots, x_n),
\end{align}
and
\begin{align}\label{e-H}
H_m(x)=\prod_{1\leqslant i < j\leqslant n}(x_i-x_j)^{2m-1}\cdot L(x_1,\ldots, x_n),
\end{align}
where $L(x_1,\ldots, x_n)$ is a symmetric polynomial in the $x$'s.

We have the following relation between the coefficients of
$H'_m(x)$ and $H_m(x)$.
\begin{lem}\label{lem-relation}
\begin{align}
[x_1^{k-1}\cdots x_n^{k-1}]H'_m(x)=n![x_1^{k-n}\cdots x_n^{k-1}]
H_m(x),
\end{align}
where $[x_1^{j_1}\cdots x_{n}^{j_n}]P(x_1,\ldots,x_n)$
denotes the coefficient of the monomial $x_1^{j_1}\cdots x_{n}^{j_n}$
in the polynomial $P(x_1,\ldots,x_n)$.
\end{lem}
\begin{proof}
By the formula of $H'_m(x)$ in \eqref{e-H'},
we have
\begin{align}
[x_1^{k-1}\cdots x_n^{k-1}]H'_m(x)=
&[x_1^{k-1}\cdots x_n^{k-1}]\prod_{1\leqslant i < j \leqslant n}(x_i-x_j)
\prod_{1\leqslant i < j \leqslant n}(x_i-x_j)^{2m-1}L(x_1,\ldots,x_n).
\end{align}
It is well-known that
$$\prod_{1\leqslant i < j \leqslant n}(x_i-x_j)=\sum_{w\in S_{n}}({\rm sgn} w)\prod_{l=1}^nx_{w(l)}^{n-l},$$
where $S_n$ is the symmetric group of all permutations on $\{1,\ldots,n\}$
and ${\rm sgn}w$ equals $1$ or $-1$ according to whether $w$ is even or odd.
Since $\prod_{1\leqslant i < j \leqslant n}(x_i-x_j)^{2m-1}L(x_1,\ldots,x_n)$
is antisymmetric in the $x$'s,
we have
{\small
\begin{align}\label{e-1}
[x_1^{k-1}\cdots x_n^{k-1}]H'_m(x)=&[x_1^{k-1}\cdots x_n^{k-1}]
\sum_{w\in S_n}({\rm sgn}w)\prod_{l=1}^nx_{w(l)}^{n-l} \nonumber \\
&\cdot({\rm sgn}w)\prod_{1\leqslant i < j \leqslant n}(x_{w(i)}-x_{w(j)})
^{2m-1}L(x_{w(1)},\ldots,x_{w(n)}) \nonumber \\
=&[x_1^{k-1}\cdots x_n^{k-1}]\sum_{w\in S_n}
\prod_{l=1}^n x_{w(l)}^{n-l}\prod_{1\leqslant i < j \leqslant n}(x_{w(i)}-x_{w(j)})^{2m-1}
L(x_{w(1)},\ldots,x_{w(n)}).
\end{align}}
Since the monomial $x_1^{k-1}\cdots x_{n}^{k-1}$
is symmetric in the $x$'s, \eqref{e-1} becomes
\begin{align*}
[x_1^{k-1}\cdots x_n^{k-1}]H'_m(x)=&n![x_1^{k-1}\cdots x_n^{k-1}] \prod_{l=1}^n x_l^{n-l}\prod_{1\leqslant i < j \leqslant n}(x_i-x_j)^{2m-1}L(x_1,\ldots,x_n) \nonumber \\
=&n![x_1^{k-n}\cdots x_n^{k-1}]H_m(x).
\end{align*}
\end{proof}
Note that this lemma can be extended to the $q$-cases, see J. Stembridge \cite[Theorem 4.1]{Stem}.

\section{Auxiliary propositions}

Denote
\begin{align}\label{def-delta}
\Delta=\prod_{l=0}^{n-1}\frac{(m(l+1))!}{(b+ml)!m!}.
\end{align}

By taking the leading and the second leading
coefficients of the Morris identity \cite{morris82}, we obtained the
following identities.
\begin{prop}\label{c-leadingM}
\emph{(Gessel et al. \cite{Zhou})}.
For $n,b,m\in \mathbb{N}$, we have
\begin{align}\label{M-leading}
\Big[\prod_{l=1}^nx_l^{b+m(n-1)}\Big]\Big(\sum_{i=1}^nx_i\Big)^{\!nb}\prod_{1\leq
i < j\leq n}(x_i-x_j)^{2m}=(-1)^{m{n\choose 2}}(nb)!\Delta,
\end{align}
and for $1\leqslant r\leqslant n$,
{\small \begin{align}\label{M-Sleading}
\Big[\prod_{l=1}^nx_l^{b+m(n-1)}\Big]&x_{r}^2\Big(\sum_{i=1}^n x_{i}\Big)^{\!nb-2}
 \prod_{1\leqslant i < j\leqslant n}(x_i-x_j)^{2m}
=(-1)^{m{n\choose 2}}b\big(b-m(n-1)-1\big)(nb-2)!\Delta.
\end{align}
}
\end{prop}
If we let $b+m(n-1)=k-1$, then \eqref{M-leading} becomes the key identity
in \cite[Proposition 2.2]{Hou-Sun} of Q. H. Hou and Z. W. Sun.
The identity led to Theorem \ref{thm-HouSun} by the polynomial method(Lemma \ref{l-polylem}).
On the other hand, by applying Lemma \ref{lem-relation} to \eqref{M-leading} we can get Corollary 2.2 in
\cite{Sun-Yeh} of Z. W. Sun and Y. N. Yeh, which led to Theorem \ref{thm-SunYeh} by Lemma \ref{l-polylem}.
Thus Theorem \ref{thm-HouSun} and Theorem \ref{thm-SunYeh} are equivalent
in the view of Lemma \ref{lem-relation}.

From \eqref{M-leading} and \eqref{M-Sleading} we can deduce the following identities.
\begin{prop}\label{prop-main}
For $1\leqslant r\leqslant n$, we have
\begin{align}\label{M-1para}
\Big[\prod_{l=1}^nx_l^{b+m(n-1)}\Big]x_r\Big(\sum_{i=1}^nx_i\Big)^{\!nb-1}\prod_{1\leq
i < j\leq n}(x_i-x_j)^{2m}=(-1)^{m{n\choose 2}}b(nb-1)!\Delta.
\end{align}
For $1\leqslant i\neq j\leqslant n$, we have
\begin{align}\label{M-2para}
\Big[\prod_{l=1}^nx_l^{b+m(n-1)}\Big]x_{i}x_{j}\Big(\sum_{i=1}^n x_{i}\Big)^{\!nb-2}
\prod_{1\leqslant i < j\leqslant n}(x_i-x_j)^{2m}
=(-1)^{m{n\choose 2}}b(b+m)(nb-2)!\Delta.
\end{align}
\end{prop}
\begin{proof}
Rewrite \eqref{M-leading} as
\begin{align*}
\Big[\prod_{l=1}^nx_l^{b+m(n-1)}\Big]\Big(\sum_{i=1}^nx_i\Big)\Big(\sum_{i=1}^nx_i\Big)^{\!nb-1}\prod_{1\leq
i < j\leq n}(x_i-x_j)^{2m}=(-1)^{m{n\choose 2}}(nb)!\Delta.
\end{align*}
Therefore by the symmetry of the $x$'s we can get \eqref{M-1para}.

Rewrite the left hand side of \eqref{M-leading} as
\begin{align*}
\Big[\prod_{l=1}^nx_l^{b+m(n-1)}\Big]\Big(\sum_{i=1}^nx_i\Big)^2 \Big(\sum_{i=1}^nx_i\Big)^{\!nb-2}\prod_{1\leq
i < j\leq n}(x_i-x_j)^{2m}.
\end{align*}
Expand $(\sum_{i=1}^nx_i)^2$ and by the symmetry of the $x$'s,
for $1\leqslant i \neq j \leqslant n$ and $1\leqslant r \leqslant n$ the above equation becomes
\begin{align*}
\Big[\prod_{l=1}^nx_l^{b+m(n-1)}\Big]\big(nx_r^2+n(n-1)x_ix_j\big)\Big(\sum_{i=1}^nx_i\Big)^{\!nb-2}\prod_{1\leq
i < j\leq n}(x_i-x_j)^{2m}.
\end{align*}
It leads to \eqref{M-2para} by substituting \eqref{M-Sleading} into the above equation.
\end{proof}

By \eqref{M-Sleading}--\eqref{M-2para} and Lemma \ref{lem-relation}
we can get the following identities.
\begin{prop}\label{prop-main2}
For $1\leqslant r\leqslant n$ and $1\leqslant i\neq j\leqslant n$, we have
{\small
\begin{align}
\Big[\prod_{l=1}^nx_l^{b+m(n-1)-n+l}\Big]x_r\Big(\sum_{i=1}^nx_i\Big)^{\!nb-1}\prod_{1\leq
i < j\leq n}(x_i-x_j)^{2m-1}=(-1)^{m{n\choose 2}}\frac{b(nb-1)!}{n!}\Delta,
\end{align}

\begin{align}
\Big[\prod_{l=1}^nx_l^{b+m(n-1)-n+l}\Big]x_{r}^2\Big(\sum_{i=1}^n x_{i}\Big)^{\!nb-2}
 \prod_{1\leqslant i < j\leqslant n}(x_i-x_j)^{2m-1}
=(-1)^{m{n\choose 2}}\frac{b\big(b-m(n-1)-1\big)(nb-2)!}{n!}\Delta,
\end{align}

\begin{align}
\Big[\prod_{l=1}^nx_l^{b+m(n-1)-n+l}\Big]x_{i}x_{j}\Big(\sum_{i=1}^n x_{i}\Big)^{\!nb-2}
\prod_{1\leqslant i < j\leqslant n}(x_i-x_j)^{2m-1}
=(-1)^{m{n\choose 2}}\frac{b(b+m)(nb-2)!}{n!}\Delta.
\end{align}
}
\end{prop}

\section{Proof of the theorems}

Our proof of Theorem \ref{thm1} follows the same line as Z. W. Sun
and Y. N .Yeh's proof of Theorem \ref{thm-SunYeh}.
\begin{proof}[Proof of Theorem {\bf \ref{thm1}}]
The case $n=1$ or $k-1<m(n-1)$ is trivial. So
we assume $n\geqslant 2$ and $b=k-1-m(n-1)\geqslant 0$.

Since $|F| \geqslant p > mn \geqslant 2m$, we can extend each
$S_{ij}(1 \leqslant i < j \leqslant n)$ to a subset $S_{ij}^{*}$
of $F$ with cardinality $2m$. By Lemma \ref{l-polylem} it suffices to show that
\begin{align*}
[x_1^{k-1}\cdots x_n^{k-1}]
\Big(\sum_{i\in T}\alpha_i x_{i}\Big)\Big(\sum_{i=1}^nx_i\Big)^{nb-1}
\prod_{1\leqslant i < j\leqslant n}\prod_{c\in S_{ij}^{*}}(x_i-x_j+c)
\end{align*}
does not vanish. Let $e$ denote the multiplicative identity
of the field $F$. Then the above coefficient equals $he$ where
$$h= [x_1^{k-1}\cdots x_n^{k-1}]
\Big(\sum_{i\in T}\alpha_i x_{i}\Big)\Big(\sum_{i=1}^nx_i\Big)^{nb-1}
\prod_{1\leqslant i < j\leqslant n}(x_i-x_j)^{2m} \in \mathbb{Z}.$$
By \eqref{M-1para},
$$h=\Big(\sum_{i\in T}\alpha_{i}\Big)(-1)^{m{n\choose 2}}b(nb-1)!\prod_{l=0}^{n-1}\frac{(m(l+1))!}{(b+ml)!m!}.$$
As $p>mn, p>nb-1$ and $\sum_{i\in T}\alpha_i\neq 0$, $p$ does not divide $h$ and hence $he\neq 0$.
This concludes the proof.
\end{proof}
We notice that recently Z. W. Sun and L. L. Zhao got a theorem \cite[Theorem 1.3]{Sun-Zhao}
to count the lower bounds for the sumsets with general polynomial restrictions,
but it seems that their theorem is not appropriate in getting Theorem \ref{thm1} directly.

The proofs of other theorems in this paper are routine
by applying Lemma \ref{l-polylem} to \eqref{M-Sleading}, \eqref{M-2para} and
Proposition \ref{prop-main2}.

\end{document}